\documentclass{article}

\usepackage{
 amsmath,
 amsxtra,
 amsthm,
 amssymb,
 etex,
 mathrsfs,
 }
\usepackage[all]{xy}

\newtheorem{theorem}{Theorem}[section]
\newtheorem{lemma}[theorem]{Lemma}

\newtheorem{proposition}[theorem]{Proposition}
\newtheorem{corollary}[theorem]{Corollary}
\newtheorem*{thm2}{Theorem}

\theoremstyle{definition}
\newtheorem{defn}[theorem]{Definition}
\newtheorem{remark}[theorem]{Remark}

\newcommand{\bd}{\begin{defn}}
\newcommand{\ed}{\end{defn}}
\newcommand{\bl}{\begin{lemma}}
\newcommand{\el}{\end{lemma}}
\newcommand{\bp}{\begin{proposition}}
\newcommand{\ep}{\end{proposition}}
\newcommand{\bt}{\begin{theorem}}
\newcommand{\et}{\end{theorem}}
\newcommand{\bc}{\begin{corollary}}
\newcommand{\ec}{\end{corollary}}
\newcommand{\br}{\begin{remark}}
\newcommand{\er}{\end{remark}}
\newcommand{\ba}{\begin{array}}
\newcommand{\ea}{\end{array}}
\newcommand{\bpf}{\begin{proof}}
\newcommand{\epf}{\end{proof}}

\newcommand{\Z}{\mathbb{Z}}
\newcommand{\Q}{\mathbb{Q}}
\newcommand{\Zp}{\mathbb{Z}_{p}}

\newcommand{\Op}{\mathcal{O}}

\newcommand{\Ap}{A[p^{\infty}]}
\newcommand{\Af}{A_f(j)}

\newcommand{\ga}{\gamma}

\newcommand{\la}{\lambda}

\newcommand{\Si}{\Sigma}

\newcommand{\Hi}{H_{\mathrm{Iw}}}

\DeclareMathOperator{\Sel}{Sel} \DeclareMathOperator{\Gal}{Gal}
\DeclareMathOperator{\Hom}{Hom} \DeclareMathOperator{\rank}{rank}
\DeclareMathOperator{\corank}{corank}
\DeclareMathOperator{\Ext}{Ext}

\newcommand{\M}{\mathfrak{M}}

\newcommand{\mK}{\mathcal{K}}

\newcommand{\ot}{\otimes}
\newcommand{\ilim}{\displaystyle \mathop{\varinjlim}\limits}
\newcommand{\plim}{\displaystyle \mathop{\varprojlim}\limits}

\newcommand{\lra}{\longrightarrow}

\newcommand{\ps}[1]{[[ #1 ]]}

  \DeclareFontFamily{U}{wncy}{}
  \DeclareFontShape{U}{wncy}{m}{n}{<->wncyr10}{}
  \DeclareSymbolFont{mcy}{U}{wncy}{m}{n}
  \DeclareMathSymbol{\sha}{\mathord}{mcy}{"58}

\begin{document}
\title{On the weak Leopoldt conjecture and coranks of Selmer groups of supersingular abelian varieties in $p$-adic Lie extensions}
 \author{Meng Fai Lim\footnote{School of Mathematics and Statistics $\&$ Hubei Key Laboratory of Mathematical Sciences, Central China Normal University, Wuhan, 430079, P.R.China. E-mail: \texttt{limmf@mail.ccnu.edu.cn}} }
\date{}
\maketitle

\begin{abstract} \footnotesize
\noindent
Let $A$ be an abelian variety defined over a number field $F$ with supersingular reduction at all primes of $F$ above $p$. We establish an equivalence between the weak Leopoldt conjecture and the expected value of the corank of the classical Selmer group of $A$ over a $p$-adic Lie extension (not neccesasily containing the cyclotomic $\Zp$-extension). As an application, we obtain the exactness of the defining sequence of the Selmer group. In the event that the $p$-adic Lie extension is one-dimensional, we show that the dual Selmer group has no nontrivial finite submodules. Finally, we show that the aforementioned conclusions carry over to the Selmer group of a non-ordinary cuspidal modular form.


\end{abstract}

\section{Introduction}

Throughout, $p$ will always denote a fixed odd prime. Let $A$ be abelian variety defined over a number field $F$ which has good supersingular reduction at all primes above $p$.  Let $F_{\infty}$ be a Galois extension of $F$ such that $G=\Gal(F_{\infty}/F)$ is a compact pro-$p$ $p$-adic Lie group with no $p$-torsion, that $F_{\infty}/F$ is unramified outside a set of finite primes and is ramified at every prime of $F$ above $p$. Let $S$ be a finite set of primes of $F$ which contains the primes above $p$, the
infinite primes, the primes at which the abelian variety $A$ has bad reduction and the
primes that are ramified in $F_{\infty}/F$. For every
algebraic (possibly infinite) extension $\mathcal{L}$ of $F$
contained in $F_S$, write $G_S(\mathcal{L}) =
\Gal(F_S/\mathcal{L})$, where $F_S$ is the
maximal algebraic extension of $F$ unramified outside $S$. One can define the classical $p$-primary Selmer group $\Sel(A/F_{\infty})$ of $A$ over $F_\infty$ which is endowed with a $\Zp\ps{G}$-module structure. We now consider the following two statements.

\medskip
$\mathbf{(A1)}$ $H^2(G_{S}(F_{\infty}), \Ap)=0$.

\medskip
$\mathbf{(A2)}$  $\corank_{\Zp\ps{G}}\Big(\Sel(A/F_{\infty})\Big) = g|F:\Q|$, where $g$ is the dimension of the abelian variety $A$.

\medskip
Statement $\mathbf{(A1)}$ is usually coined as the weak Leopoldt conjecture. When the extension $F_{\infty}$ contains the cyclotomic $\Zp$-extension, statement $\mathbf{(A2)}$ is a special case of a conjecture of Schneider \cite{Sch85} which is a generalization of a conjecture of Mazur \cite{Maz}. In this note, we shall show that the two statements are equivalent.

\begin{thm2}[Theorem \ref{WLSch}]
 Let $A$ be a $g$-dimensional abelian variety defined over a number field $F$ which has good supersingular reduction at all primes above $p$. Let $F_{\infty}$ be a Galois extension of $F$ such that $\Gal(F_{\infty}/F)$ is a compact pro-$p$ $p$-adic Lie group with no $p$-torsion, that $F_{\infty}/F$ is unramified outside a set of finite primes and is ramified at every prime of $F$ above $p$.

 Then $H^2(G_S(F_{\infty}),\Ap)=0$ if and only if
\[\corank_{\Zp\ps{G}}\Big(\Sel(A/F_{\infty})\Big) = g|F:\Q|.\]
\end{thm2}

In the situation that $F_\infty$ contains the cyclotomic $\Zp$-extension, the implication
$\mathbf{(A1)}\Rightarrow \mathbf{(A2)}$ has been observed in \cite{C99} and somewhat implicitly in \cite{CM, Kob, Mc}. When $F_\infty$ is the anti-cyclotomic $\Zp$ of an imaginary field, such equivalence is also examined in \cite{Ci, LV}. One point we like to stress is that our extension $F_\infty$ is not assumed to contain the cyclotomic $\Zp$-extension. This does not correlate to the situation when the abelian variety has some primes above $p$ which is not a supersingular prime. For instance, if we consider an elliptic curve $E$ defined over an imaginary quadratic field $K$ that has \textit{good ordinary reduction} at all primes above $p$, a conjecture of Mazur \cite{Maz} predicts that the Selmer group of $E$ over the cyclotomic $\Zp$-extension of $K$ is cotorsion (this is a theorem of Kato \cite{K} if $E$ is defined over $\Q$). On the other hand, if one considers the Selmer group of $E$ over the anticyclotomic $\Zp$-extension of $K$, the said Selmer group can be \textit{non-cotorsion} (for instance, see \cite{Be, Maz84}). Our theorem is essentially saying that such disparity of the corank of the Selmer group of a supersingular abelian variety in different classes of $p$-adic Lie extensions will not happen (modulo weak Leopoldt).

Another point we like to mention is that when some of the prime above $p$ is not a supersingular prime of $A$, one usually has an assertion saying corank of the Selmer group over an $p$-adic Lie extension takes a certain predicted value if and only if the defining sequence for the Selmer group is short exact and the weak Leopoldt conjecture is valid (for instance, see \cite[Theorem 4.12]{C99} or \cite[Proposition 3.3]{LimMHG}). Pertaining to this aspect, our result is saying that for supersingular abelian variety, the short exactness of the defining sequence does not come into play.

We now give an outline of the paper.
The above mentioned result will be proved in Section \ref{main section}. We shall also give various sufficient criterions to verify the validity of  $\mathbf{(A1)}$ and $\mathbf{(A2)}$ (see Propositions \ref{WLSch coates} and \ref{WLSch solvable}, and Corollary \ref{WLSch cor}).
In Section \ref{surjective section}, we return to study the question of the exactness of the defining sequence of the Selmer group. Namely, we show that this is a consequence of $\mathbf{(A1)}$ (see Proposition \ref{surjective}). This again seems to be a phenomenon reserved for supersingular abelian varieties. We then apply this exactness property to show that the dual Selmer group has no non-zero pseudo-null submodules (see Theorem \ref{no finite submodule theorem}) in Section \ref{non-existence section}. This is an exact analog of a result of Greenberg \cite[Proposition 4.14]{G99}, where he obtained the same conclusion for the Selmer group of an elliptic curve with good ordinary reduction at all primes above $p$. We also note that our result refines \cite[Theorem 2.14]{LeiSu}, where they establish a similar conclusion with slightly more stringent hypothesis.
Finally, in Section \ref{non-ordinary section}, we consider the Selmer group of a non-ordinary modular form over a $\Zp$-extension and show that the aforementioned results carry over to this context.

\subsection*{Acknowledgement}
The author would like to thank Antonio Lei for
many insightful discussions and for
answering many questions on his papers \cite{Lei11, LeiSu} . This research is supported by the National Natural Science Foundation of China under Grant No. 11550110172 and Grant No. 11771164.

\section{Main result} \label{main section}

Let $A$ be a $g$-dimensional abelian variety defined over a number field $F$ which has supersingular reduction at all primes above $p$. Let $v$ be a prime of $F$. For every finite
extension $L$ of $F$, we define
 \[ J_v(A/L) = \bigoplus_{w|v}H^1(L_w, A)[p^\infty],\]
where $w$ runs over the (finite) set of primes of $L$ above $v$. If
$\mathcal{L}$ is an infinite extension of $F$, we define
\[ J_v(A/\mathcal{L}) = \ilim_L J_v(A/L),\]
where the direct limit is taken over all finite extensions $L$ of
$F$ contained in $\mathcal{L}$. For any algebraic (possibly
infinite) extension $\mathcal{L}$ of $F$, the $p$-primary Selmer group of $A$
over $\mathcal{L}$ is defined to be
\[ \Sel(A/\mathcal{L}) = \ker\Big(H^1(\mathcal{L}, \Ap)\lra \bigoplus_{v} J_v(A/\mathcal{L})
\Big), \]
 where $v$ runs through all the primes of $F$.

Let $S$ be a finite set of primes of $F$ which contains the primes above $p$, the
infinite primes, the primes at which the abelian variety $A$ has bad reduction and the
primes that are ramified in $F_{\infty}/F$. Denote by $F_S$ the
maximal algebraic extension of $F$ unramified outside $S$. For every
algebraic (possibly infinite) extension $\mathcal{L}$ of $F$
contained in $F_S$, write $G_S(\mathcal{L}) =
\Gal(F_S/\mathcal{L})$
. The following alternative equivalent
description of the Selmer group
\[ \Sel(A/\mathcal{L}) = \ker\Big(H^1(G_S(\mathcal{L}), \Ap)\stackrel{\la_S(F_{\infty})}{\lra} \bigoplus_{v\in S} J_v(A/\mathcal{L})
\Big)\]
is well-known (see \cite[Chap. I, Corollary 6.6]{Mi}). Denote by $X(A/\mathcal{L})$ the
Pontryagin dual of $\Sel(A/\mathcal{L})$.
We now introduce the extension which we are interested in.

\bd A Galois extension $F_{\infty}$ of $F$ is said to be a strongly pro-$p$ $p$-adic Lie
extension of $F$ if all of the following statements hold.

(1) $G:=\Gal(F_{\infty}/F)$ is a compact pro-$p$ $p$-adic
Lie group with no $p$-torsion.

(2) The extension $F_{\infty}$ is unramified outside a
finite set of primes of $F$.

(3) All the primes of $F$ above $p$ are ramified in $F_{\infty}/F$.
\ed

Before continuing our discussion, we quickly review the notion of $\Zp\ps{G}$-rank.
Since the group $G=\Gal(F_{\infty}/F)$ is a compact pro-$p$ $p$-adic
Lie group without $p$-torsion, it follows that $\Zp\ps{G}$ is
an Auslander regular ring (cf. \cite[Theorems 3.26]{V02}) and has no zero divisors (cf.\
\cite{Neu}). Consequently, the ring $\Zp\ps{G}$ admits a skew field $Q(G)$ which is flat
over $\Zp\ps{G}$ (see \cite[Chapters 6 and 10]{GW} or \cite[Chapter
4, \S 9 and \S 10]{Lam}). This allows us to define the notion of $\Zp\ps{G}$-rank of a finitely generated $\Zp\ps{G}$-module $M$, which is given by
\[\rank_{\Zp\ps{G}}(M)  = \dim_{Q(G)} (Q(G)\ot_{\Zp\ps{G}}M). \]
If $N$ is a discrete $\Zp\ps{G}$-module such that its Pontryagin dual $N^\vee$ is finitely generated over $\Zp\ps{G}$, we write
\[\corank_{\Zp\ps{G}}(N)  = \rank_{\Zp\ps{G}}(N^{\vee}). \]

We turn back to our arithmetic situation. Consider the following two statements as stated in the introduction.

\medskip
$\mathbf{(A1)}$ $H^2\big(G_{S}(F_{\infty}), \Ap\big)=0$.

\medskip
$\mathbf{(A2)}$  $\corank_{\Zp\ps{G}}\Big(\Sel(A/F_{\infty})\Big) = g|F:\Q|$, where $g$ is the dimension of $A$.

\medskip
As mentioned in the introduction, the first statement is usually called the weak Leopoldt conjecture. When $F_{\infty}$ contains the cyclotomic $\Zp$-extension of $F$, statement $\mathbf{(A2)}$ is a special case of a conjecture of Schneider \cite{Sch85} (also see \cite{C99, OcV03}). In this note, we shall work with general extensions which may not contain the cyclotomic $\Zp$-extension.
The following is the theorem as stated in the introduction.

\bt \label{WLSch}
 Let $A$ be a $g$-dimensional abelian variety defined over a number field $F$ which has supersingular reduction at all primes above $p$. Let $F_{\infty}$ be a strongly pro-$p$ $p$-adic extension of $F$.

 Then we have $H^2(G_S(F_{\infty}),\Ap)=0$ if and only if
\[\corank_{\Zp\ps{G}}\Big(\Sel(A/F_{\infty})\Big) = g|F:\Q|.\]
\et

Before proving the theorem, we require a preparatory lemma.

\bl \label{WLSch lemma}
For every $v|p$, we have $J_v(A/F_{\infty})=0$. For every $v\nmid p$, $J_v(A/F_{\infty})$ is cotorsion over $\Zp\ps{G}$.
\el

\bpf
For each $v\in S$, write $w$ for a prime of $F_\infty$ above $v$. Now suppose that $v|p$. By \cite[P. 150]{CG}, we have the following short exact sequence
\[ 0\lra C_v \lra \Ap\lra D_v\lra 0\]
of $\Gal(\bar{F_v}/F_v)$-modules, where $C_v$ and $D_v$ are cofree $\Zp$-modules  characterized by the facts that $C_v$ is divisible and that $D_v$ is the maximal quotient of $\Ap$ by a divisible subgroup such that $\Gal(\bar{F_v}/F_v^{ur})$ acts on $D_v$ via a finite quotient. By the hypothesis on $F_\infty$, the prime $v$ is infinitely ramified in $F_\infty$ and so $F_{\infty,w}/F_v$ is infinitely ramified. Therefore, we may apply \cite[Proposition 4.8]{CG} to conclude that
$H^1(F_{\infty,w},A)[p^{\infty}]\cong H^1(F_{\infty,w},D_v)$. Since the abelian variety $A$ is assumed to have supersingular reduction at $v$, we have $D_v=0$ and so $H^1(F_{\infty,w},A)[p^{\infty}]=0$. This in turn implies that $J_v(A/F_{\infty})=0$.

Now, suppose that $v$ does not divide $p$. We first consider the situation that $v$ does not decompose completely in $F_\infty/F$. Write $G_w$ for the decomposition group of $w$ in $G$ which is nontrivial by our assumption. Since $G$ has no $p$-torsion, it follows that the dimension of $G_w$ is at least one. Thus, we have
\[J_v(A/F_\infty) = \mathrm{Coind}_{G_w}^{G}\Big(H^1(F_{\infty,w},\Ap)\Big).\]
By \cite[Theorem 4.1]{OcV03}, $H^1(F_{\infty,w},\Ap)$ is cotorsion over $\Zp\ps{G_w}$ and hence $J_v(A/F_\infty)$ is cotorsion over $\Zp\ps{G}$. We now suppose that the prime $v$ decomposes completely in $F_\infty/F$. In this case, one has $H^1(F_{\infty,w},\Ap) = H^1(F_v,\Ap)$, By the Tate duality (cf. \cite[Chap.\ I, Corollary 3.4]{Mi}), the latter is isomorphic to $\plim_m A^t(F_v)/p^m$, where $A^t$ is the dual abelian variety of $A$. Now, a well-known theorem of Mattuck \cite{Mat} tells us that the group $A^t(F_v)$ is finitely
generated over $\Z_l$, where $l\neq p$. It follows from this that $\plim_m A^t(F_v)/p^m$ is finite. In conclusion, we see that for large enough integer $t$, $p^t$ annihilates $H^1(F_{\infty,w},\Ap)^\vee$ and hence $J_v(A/F_{\infty})^{\vee}$. Thus,  $J_v(A/F_\infty)$ is also cotorsion over $\Zp\ps{G}$ in this case. This completes the proof of the lemma.
\epf

We can now give the proof of Theorem \ref{WLSch}.

\bpf[Proof of Theorem \ref{WLSch}]
 By virtue of Lemma \ref{WLSch lemma} and the finiteness of $S$, we see that
  \[\bigoplus_{v\in S} J_v(A/F_\infty)\]
 is cotorsion over $\Zp\ps{G}$. It then follows from this that
 \[\corank_{\Zp\ps{G}}\Big(\Sel(A/F_{\infty})\Big) = \corank_{\Zp\ps{G}}\Big(H^1(G_{\Si}(F_{\infty}),\Ap)\Big). \]
On the other hand, a standard Euler-characteristics argument (cf.\ \cite[Proposition 3]{G89} or \cite[Theorem 3.2]{OcV03}) yields
 \[ \corank_{\Zp\ps{G}}\Big(H^1(G_{\Si}(F_{\infty}),\Ap)\Big)- \corank_{\Zp\ps{G}}\Big(H^2(G_{\Si}(F_{\infty}),\Ap)\Big)  = g|F:\Q|. \]
Combining these calculations, we have that
 \[\rank_{\Zp\ps{G}}\Big(X(A/F_{\infty})\Big) = g|F:\Q|\]
 if and only if one has
  \[ \corank_{\Zp\ps{G}}\Big(H^2(G_{\Si}(F_{\infty}),\Ap)\Big)=0.\]
  On the other hand, it is well-known that $H^2(G_{\Si}(F_{\infty}),\Ap)^{\vee}$ is a submodule of a projective $\Zp\ps{G}$-module (cf.\ \cite[Diagram 2.1]{OcV03}) and hence must be torsionfree over $\Zp\ps{G}$. Hence the latter corank equality is precisely equivalent to $H^2(G_{S}(F_{\infty}),\Ap)=0$. This completes the proof of the theorem.
\epf

\br The proof is inspired by \cite[Corollary 1.9]{C99}. To the best of the author's knowledge, the implication of the theorem does not hold in the event that $E$ has non-supersingular reduction at some prime above $p$. One usually needs an additional input on the surjectivity of $\rho$ (for instance, see \cite[Theorem 4.12]{C99} or \cite[Proposition 3.3]{LimMHG}).
\er

For the remainder of this section, we discuss some situations, where one can establish the validity of $\mathbf{(A1)}$ and $\mathbf{(A2)}$. As a start, we have the following observation which goes back to \cite[Corollary 1.9]{C99}.

\bp \label{WLSch coates}
Let $A$ be a $g$-dimensional abelian variety defined over a number field $F$ which has supersingular reduction at all primes above $p$. Let $F_{\infty}$ be a strongly pro-$p$ $p$-adic extension of $F$ such that $\Gal(F_S/F_\infty)$ acts trivially on $\Ap$.

 Then both assertions $\mathbf{(A1)}$ and $\mathbf{(A2)}$ hold.
\ep

\bpf
 It follows from an observation of Ochi (see \cite[Theorem 2.10]{C99} or \cite[Proof of Corollary 4.8]{OcV03}) that under the hypothesis that $\Gal(F_S/F_\infty)$ acts trivially on $\Ap$, one has $H^2(G_{S}(F_{\infty}),\Ap)=0$. The conclusion of the proposition now follows from Theorem \ref{WLSch}.
\epf

We give another criterion for the validity of $\mathbf{(A1)}$ and $\mathbf{(A2)}$.

\bp \label{WLSch solvable}
 Let $A$ be a $g$-dimensional abelian variety defined over a number field $F$ which has good supersingular reduction at all primes above $p$. Let $F_{\infty}$ be a strongly pro-$p$ $p$-adic extension of $F$ such that $\Gal(F_{\infty}/F)$ is solvable. Suppose that $H^2(G_S(F),\Ap)=0$.

 Then we have $H^2(G_S(F_{\infty}),\Ap)=0$ and
\[\corank_{\Zp\ps{G}}\Big(\Sel(A/F_{\infty})\Big) = g|F:\Q|.\]
\ep

Write $T_pA$ for the Tate module of the abelian variety $A$. Define $\Hi^i(F_\infty/F, T_pA) = \plim_L H^i(G_S(L), T_pA)$, where the inverse limit is taken over all finite subextensions of $F_{\infty}/F$ and the transition maps are given by the corestriction maps. We now prove the following lemma which will be required in the proof of Proposition \ref{WLSch solvable}.

\bl \label{Jannsen}
Let $A$ be a $g$-dimensional abelian variety defined over a number field $F$ which has good supersingular reduction at all primes above $p$. Let $F_{\infty}$ be a strongly pro-$p$ $p$-adic extension of $F$. Then $H^2(G_S(F_\infty), \Ap)=0$ if and only if $\Hi^2(F_\infty/F, T_pA)$ is torsion over $\Zp\ps{G}$.
\el

\bpf
This follows from a standard argument appealing to the spectral sequence of Jannsen (for instance, see \cite[Theorem 3.3]{OcV03} or \cite[Lemma 7.1]{LimFine}). For the convenience of the readers, we supply the argument here.
 Consider the spectral sequence of Jannsen
\[ E^{i,j}_2=\Ext^i_{\Zp\ps{G}}\big(H^j(G_S(F_\infty),A[p^\infty])^{\vee},\Zp\ps{G}\big) \Longrightarrow \Hi^{i+j}(F_{\infty}/F, T_pA)\]
(cf.\ \cite[Theorem 1]{Jan}).
As the ring $\Zp\ps{G}$ has finite projective dimension, the spectral sequence is bounded. Therefore, the terms $E_m^{i,j}$ must stabilize for large enough $m$. In particular, we have that $E_{\infty}^{i,j}$
    is a subquotient of $E_2^{i,j}$. By \cite[Proposition 3.5(iii)(a)]{V02}, the terms $E_2^{i,j}$, and hence $E_{\infty}^{i,j}$, are torsion $\Zp\ps{G}$-modules for $i\neq
    0$. Since $\Hi^{2}(F_{\infty}/F, T_pA)$ has a finite filtration with
    factors $ E_{\infty}^{i,2-i}$ for $i =0,1,2$, this in turn implies that
     $ \rank_{\Zp\ps{G}}\Hi^2(F_{\infty}/F, T_pA) = \rank_{\Zp\ps{G}}E_{\infty}^{0,2}.$
    On the other hand, one sees easily that the edge map
    $E_{\infty}^{0,2}\to E_2^{0,2}$ is injective with cokernel isomorphic to a subquotient of
     \[
        \Ext_{\Zp\ps{G}}^{2}\big(H^1(G_S(F_\infty),\Ap)^{\vee}, \Zp\ps{G}\big)\oplus \Ext_{\Zp\ps{G}}^{3}\big(H^0(G_S(F_\infty),\Ap)^{\vee},
        \Zp\ps{G}\big).
    \]
    Again, it follows from \cite[Proposition 3.5(iii)(a)]{V02} that these are torsion over $\Zp\ps{G}$. Therefore, we may
    conclude that
      \[\rank_{\Zp\ps{G}}\Hi^2(F_{\infty}/F, T_pA) = \rank_{R\ps{G}}\Hom_{\Zp\ps{G}}\big( H^2(G_S(F_{\infty}),\Ap)^{\vee},
  \Zp\ps{G}\big). \]
  It then follows that $\Hi^2(F_{\infty}/F, T_pA)$ is a torsion $\Zp\ps{G}$-module
  if and only if
  \[\Hom_{\Zp\ps{G}}\Big(H^2(G_S(F_{\infty}),\Ap)^{\vee},
  \Zp\ps{G}\Big) \]
  is a torsion $\Zp\ps{G}$-module. As seen in the proof of Theorem \ref{WLSch}, $H^2(G_S(F_{\infty}),\Ap)^{\vee}$ is $\Zp\ps{G}$-torsionfree. Hence
  the latter statement holds if and only if $H^2(G_S(F_{\infty}),\Ap) = 0$. This proves the lemma.
\epf

\br
By a similar argument to that in Lemma \ref{Jannsen}, one can also show that  $H^2(G_S(F),\Ap)=0$ if and only if $H^2(G_S(F),T_pA)$ is finite.
\er

We can now prove Proposition \ref{WLSch solvable}.

\bpf[Proof of Proposition \ref{WLSch solvable}]
By the hypothesis that $H^2(G_S(F),\Ap)=0$ and the above remark, we have that $H^2(G_S(F),T_pA)$ is finite. By considering the initial term of the homological spectral
sequence
\[ H_{i}\Big(G, \Hi^{-j}(F_{\infty}/F,T_pA)\Big) \Longrightarrow H^{-i-j}\big(G_S(F), T_pA\big) \]
(cf. \cite[Theorem 3.1.8]{LSh}), we have
$\Hi^2(F_\infty/F,T_pA)_G\cong H^2(G_S(F),T_pA)$. Since $G$ is solvable, it follows from \cite[pp 229 Theorem]{BH} and the above isomorphism that $\Hi^2(F_\infty/F,T_pA)$ is torsion over $\Zp\ps{G}$. By Lemma \ref{Jannsen}, the latter is equivalent to $H^2(G_S(F_\infty),\Ap)=0$. The remaining conclusion of the proposition then follows from Theorem \ref{WLSch}.
\epf

We end the section recording the following corollary.

\bc \label{WLSch cor}
Retain the assumptions of the preceding theorem. Suppose  further that the following assertions hold.

$(a)$ $\Gal(F_{\infty}/F)$ is solvable.

$(b)$ The $p$-primary part of the Tate-Shafarevich group $\sha(A/F)[p^{\infty}]$ is finite.

$(c)$ $\rank_{\Z}(A(F))\leq 1$.

Then both statements $\mathbf{(A1)}$ and $\mathbf{(A2)}$ are valid.
\ec

\bpf
By a similar argument to that in \cite[Theorem 12]{CM}, one can show that $H^2(G_S(F),\Ap)=0$ under the validity of (b) and (c). The conclusion then follows from an application of Proposition \ref{WLSch solvable}.
\epf

\section{Surjectivity of the localisation map} \label{surjective section}

Retain the notation of the previous section. We will study the surjectivity of the localization map in the definition of the Selmer groups. In subsequent discussion, we write $A^t$ for the dual abelian variety of $A$. Our result is as follow.

\bp \label{surjective}
 Let $A$ be a $g$-dimensional abelian variety defined over a number field $F$ which has good supersingular reduction at all primes above $p$. Let $F_{\infty}$ be a strongly pro-$p$ $p$-adic extension of $F$ such that $H^2(G_{S}(F_\infty),\Ap)=0$. In the event that $\Gal(F_{\infty}/F)\cong\Zp$, suppose further that $A^t(F_\infty)[p^\infty]$ is finite.
  Then we have a short exact sequence
 \[ 0\lra \Sel(A/F_\infty)\lra H^1\big(G_S(F_{\infty}),\Ap\big)\stackrel{\rho}{\lra} \bigoplus_{v\in S}J_v(A/F_\infty) \lra 0.\]
\ep

\bpf
  By \cite[Proposition A.3.2]{PR00}, we have an exact sequence
 \[ 0\lra \Sel(A/F_\infty)\lra H^1\big(G_S(F_{\infty}),\Ap\big)\stackrel{\rho}{\lra} \bigoplus_{v\in S}J_v(A/F_\infty)\lra \mathfrak{S}(A^t/F_{\infty})^{\vee} \lra 0,\]
where $\mathfrak{S}(A^t/F_\infty)$ is a $\Zp\ps{G}$-submodule of $\Hi^1(F_{\infty}/F, T_pA^t)$ and the final zero follows from the hypothesis that $H^2(G_S(F_{\infty}),\Ap)=0$.  By Lemma \ref{WLSch lemma}, $\displaystyle\bigoplus_{v\in S}J_v(A/F_\infty)$ is cotorsion over $\Zp\ps{G}$. It follows from this and the above exact sequence that $\mathfrak{S}(A^t/F_{\infty})$ is torsion over $\Zp\ps{G}$.
On the other hand, by considering the low degree terms of the spectral sequence of Jannsen
\[ \Ext^i_{\Zp\ps{G}}\big(H^j(G_S(F_\infty),A^t[p^\infty])^{\vee},\Zp\ps{G}\big) \Longrightarrow \Hi^{i+j}(F_{\infty}/F, T_pA^t)\]
(cf. \cite[Theorem 1]{Jan}), we obtain the following exact sequence
\[ 0\lra \Ext^1_{\Zp\ps{G}}\big(A^t(F^{\infty})[p^\infty])^{\vee},\Zp\ps{G}\big) \lra \Hi^1(F_{\infty}/F, T_pA^t) \lra \Ext^0_{\Zp\ps{G}}\big(H^1(G_S(F_{\infty}),A^t[p^\infty])^{\vee},\Zp\ps{G}\big). \]
Now, if $\dim G\geq 2$ and noting that $A^t(F^{\infty})[p^\infty])^{\vee}$ is finitely generated over $\Zp$, it follows from \cite[Proposition 2.2]{OcV03} that $\Ext^1_{\Zp\ps{G}}\big(A^t(F^{\infty})[p^\infty])^{\vee},\Zp\ps{G}\big)=0$. In the event that $\dim G=1$, it follows from the extra hypothesis of the proposition that  $A^t(F^{\infty})[p^\infty])^{\vee}$ is finite and so is pseudo-null over $\Zp\ps{G}$. Hence either way, the leftmost term vanishes which in turn implies that $\Hi^1(F_{\infty}/F, T_pA^t)$ injects into an $\Ext^0$-term. Since this latter term is a reflexive $\Zp\ps{G}$-module by \cite[Proposition 3.11]{V02}, it follows that $\Hi^1(F_{\infty}/F, T_pA^t)$ must be $\Zp\ps{G}$-torsionfree. Consequently, so is $\mathfrak{S}(A^t/F_{\infty})$. But as seen above, $\mathfrak{S}(A^t/F_{\infty})$ is also torsion over $\Zp\ps{G}$. Hence we must have $\mathfrak{S}(A^t/F_{\infty})=0$ and so the map $\rho$ is surjective as required.
\epf

We record a corollary of the proposition.

\bc \label{surjective cor}
 Let $A$ be a $g$-dimensional abelian variety defined over a number field $F$ which has good supersingular reduction at all primes above $p$. Let $F_{\infty}$ be a strongly pro-$p$ $p$-adic extension of $F$ such that $G=\Gal(F_{\infty}/F)$ is solvable. In the event that $\Gal(F_{\infty}/F)\cong\Zp$, suppose further that $A^t(F_\infty)[p^\infty]$ is finite.
 Assume that $H^2(G_{S}(F),\Ap)=0$.  Then we have a short exact sequence
 \[ 0\lra \Sel(A/F_{\infty})\lra H^1(G_S(F_{\infty}),\Ap)\stackrel{\rho}{\lra} \bigoplus_{v\in S}J_v(A/F_\infty)\lra 0.\]
\ec

\bpf
This follows from a combination of Propositions \ref{WLSch solvable} and \ref{surjective}.
\epf

We end with the following remark.

\br
We say a bit on the extra condition of $A^t(F_\infty)[p^\infty]$ being finite when $\dim G=1$. In the case when $F_\infty$ is the cyclotomic $\Zp$-extension, the finiteness is well-known (see \cite{Im, Ri}). For an arbitrary $\Zp$-extension, Wingberg has worked out precisely when this finiteness property holds for a simple abelian variety (see \cite[Theorem 4.3]{Win87}).
For a general abelian variety, he has also established that the finiteness fails for only finitely many $\Zp$-extensions (see \cite[Theorem 3.5]{Win87}). In view of these, the extra finiteness assumption imposed when $\dim G=1$ seems mild.
\er

\section{Non-existence of finite submodules} \label{non-existence section}

Retain the notations and assumptions from the previous sections.
In this section, we prove an analogous result of Greenberg \cite[Proposition 4.14]{G99}, where he obtained the conclusion for the Selmer group of an elliptic curve with good ordinary reduction at all primes above $p$. We also note that our result refines \cite[Theorem 2.14]{LeiSu}, where they also obtained the same conclusion for an elliptic curve with good supersingular reduction at all primes above $p$.

\bt \label{no finite submodule theorem}
 Let $A$ be a $g$-dimensional abelian variety defined over a number field $F$ which has good supersingular reduction at all primes above $p$. Let $F_{\infty}$ be a strongly pro-$p$ $p$-adic extension of $F$ of dimension one. Suppose that $A^t(F_\infty)[p^\infty]=0$ and that $H^2(G_{S}(F_\infty),\Ap)=0$.
 Then $X(A/F_\infty)$ has no non-trivial finite $\Zp\ps{G}$-submodule.
\et

\br
For $\dim G\geq 2$, the question of $X(A/F_\infty)$ having no non-trivial pseudo-null  $\Zp\ps{G}$-submodule is not so interesting (see \cite[Line after Theorem 5.1]{OcV02}). In fact, in this situation, one usually has $\Sel(A/F_\infty) = H^1(G_S(F_\infty),\Ap)$ and it is well-known that $H^1(G_S(F_\infty),\Ap)^\vee$ has no non-trivial pseudo-null  $\Zp\ps{G}$-submodule under the validity of weak Leopoldt (see \cite[Theorem 4.7]{OcV02}). Therefore, we only focus on the dimension $1$ case in this article.
\er

We first prove a special case of Theorem \ref{no finite submodule theorem}.

\bp \label{no finite submodule}
 Let $A$ be a $g$-dimensional abelian variety defined over a number field $F$ which has good supersingular reduction at all primes above $p$. Let $F_{\infty}$ be a strongly pro-$p$ $p$-adic extension of $F$ with $\Gal(F_\infty/F)\cong\Zp$. Suppose that $A^t(F_\infty)[p^\infty]=0$ and that $H^2(G_{S}(F),\Ap)=0$.
 Then $X(A/F_\infty)$ has no non-trivial finite $\Zp\ps{G}$-submodule.
\ep

We emphasis the difference lies in that we are assuming $H^2(G_{S}(F),\Ap)=0$ in Proposition \ref{no finite submodule} which is a stronger assumption that $H^2(G_{S}(F_\infty),\Ap)=0$ as assumed in Theorem \ref{no finite submodule theorem}. As a start, we establish the following lemma.

\bl \label{no finite submodule lemma}
Retain the assumptions of Proposition \ref{no finite submodule}. Then the following assertions are valid.

$(a)$ $H^1\big(G_S(F), T_pA^t\big)$ is $\Zp$-torsionfree.

$(b)$ $H^1\Big(G, H^1(G_S(F_\infty),\Ap)\Big)=0$
\el

\bpf
By considering the long cohomological exact sequence of the following short exact sequence
\[ 0\lra T_pA^t \stackrel{\cdot p}{\lra} T_pA^t\lra A^t[p]\lra 0,\]
we obtain
\[ A^t(F)[p]\lra H^1(G_S(F), T_pA^t) \stackrel{\cdot p}{\lra} H^1(G_S(F), T_pA^t). \]
Since $A^t(F)[p]=0$ by hypothesis, we have that $H^1(G_S(F), T_pA^t)$ has no $p$-torsion and this proves assertion (a).

Taking $\dim G=1$ into account, the spectral sequence
\[ H^i\Big(G, H^j(G_S(F_\infty),\Ap)\Big)\Longrightarrow H^{i+j}(G_S(F),\Ap)\]
yields an isomorphism
\[ H^1\Big(G, H^1(G_S(F_\infty),\Ap)\Big)\cong H^2(G_S(F),\Ap)\]
where the latter vanishes by hypothesis. Hence this proves assertion (b).
\epf

We can now prove Proposition \ref{no finite submodule}.

\bpf[Proof of Proposition \ref{no finite submodule}]
By \cite[Proposition A.3.2]{PR00}, we have an exact sequence
 \[ 0\lra \Sel(A/F)\lra H^1\big(G_S(F),\Ap\big)\lra \bigoplus_{v\in S}J_v(A/F)\lra \mathfrak{S}(A^t/F)^{\vee} \lra 0,\]
where $\mathfrak{S}(A^t/F)$ is a $\Zp$-submodule of $H^1(G_S(F), T_pA^t)$ and the final zero follows from the hypothesis that $H^2(G_S(F),\Ap)=0$.
On the other hand, it follows from Corollary \ref{surjective cor} that there is a short exact sequence
\[ 0\lra \Sel(A/F_{\infty})\lra H^1(G_S(F_{\infty}),\Ap)\lra \bigoplus_{v\in S}J_v(A/F_\infty)\lra 0.\]
Taking $G$-invariant and taking Lemma \ref{no finite submodule lemma}(b) into account, we have the following exact sequence
\[ 0\lra \Sel(A/F_{\infty})^G\lra H^1(G_S(F_{\infty}),\Ap)^G\lra \Big(\bigoplus_{v\in S}J_v(A/F_\infty)\Big)^G\lra H^1\Big(G,\Sel(A/F_\infty)\Big)\lra 0.\]
All of these fit into the following commutative diagram
\[   \entrymodifiers={!! <0pt, .8ex>+} \SelectTips{eu}{}\xymatrix{
       H^1(G_S(F),\Ap)
    \ar[d] \ar[r] & \displaystyle\bigoplus_{v\in S}J_v(A/F) \ar[d]_{g}  \ar[r] & \mathfrak{S}(A^t/F)^\vee \ar[d]_h  \ar[r] & 0\\
     H^1(G_S(F_{\infty}),\Ap)^G \ar[r] & \left(\displaystyle\bigoplus_{v\in S}J_v(A/F_\infty)\right)^G
    \ar[r] & H^1\Big(G,\Sel(A/F_\infty)\Big)\ar[r] &0 }\]
 with exact rows, where the leftmost two vertical maps are the restriction maps on cohomology and the rightmost map $h$ is in turn induced by these restriction maps. Since $\dim G=1$, the map $g$ is surjective and hence so is the map $h$. But as $\mathfrak{S}(A^t/F)$ is a $\Zp$-submodule of $H^1(G_S(F), T_pA^t)$, it follows from Lemma \ref{no finite submodule lemma}(a) that $\mathfrak{S}(A^t/F)^{\vee}$ is $p$-divisible. It then follows from this and the surjectivity of $h$ that $H^1\Big(G,\Sel(A/F_\infty)\Big)$ is also $p$-divisible. Therefore, we may apply \cite[Proposition 5.3.19(i)]{NSW} to obtain the required conclusion of the proposition.
\epf

We now come to the proof of Theorem \ref{no finite submodule theorem}. Before continuing, we recall the definition of the fine Selmer group over $F$ (cf. \cite{CS05, LimFine})

\[R(A/F) = \ker\Big(H^1(G_S(F), \Ap)\lra \bigoplus_{v\in S}H^1(F_v, \Ap)\Big).\]

One has a similar definition over the intermediate extensions of $F_\infty/F$. The fine Selmer group $R(A/F_\infty)$ is then defined to be the direct limit of these intermediate fine Selmer groups.  Fix an isomorphism $\kappa:G\cong 1+p\Zp$. For each $s$, we write $\Zp(s)$ for the abelian group $\Zp$ with a $G$-action given by $\ga\cdot x = \kappa(\ga)x$ for $\ga\in G$ and $x\in \Zp$. For a $G$-module $M$, write $M(s) = M\ot \Zp(s)$, where $G$ acts diagonally. To simplify notation, we shall write $A_s : = \Ap(s)$. We define $R(A_s/F)$ and $R(A_s/F_\infty)$  by replacing $\Ap$ by $A_s$ in the definition of the fine Selmer groups. Since $G_S(F_\infty)$ acts trivially on $\Zp(s)$, we have
$R(A_s/F_\infty) = R(A/F_\infty)(s)$.

The Selmer group of $A_s$ over $F_\infty$ (denoted by $\Sel(A_s/F_\infty)$) is defined similarly, with the only slight difference in that we set $J_v(A_s/F_\infty)=0$ for primes above $p$. Taking Lemma \ref{WLSch lemma} into account, we also have $\Sel(A_s/F_\infty) = \Sel(A/F_\infty)(s)$.

\bl \label{twist}
Retain the assumptions of Theorem \ref{no finite submodule theorem}.
 Then there exists $s\in \Zp$ such that $H^2(G_S(F),A_s)=0$.
\el

\bpf
By Lemma \ref{Jannsen} and the hypothesis that $H^2(G_S(F_\infty),\Ap)=0$, it follows that $\Hi^2(F_\infty/F, T_pA^t)$ is torsion over $\Zp\ps{G}$. Since $R(A/F_\infty)^{\vee}$ is contained in $\Hi^2(F_\infty/F, T_pA^t)$ by the Poitou-Tate sequence, it follows that $R(A/F_\infty)$ is also cotorsion over $\Zp\ps{G}$. As seen in the proof of \cite[Proposition 4.14]{G99}, we can find $s\in \Zp$ such that
$R(A_s/F_\infty)^G = \Big(R(A/F_\infty)(s)\Big)^G$ is finite. A standard argument then shows that the kernel of the map
\[r:R(A_s/F) \lra R(A_s/F_\infty)^G\]
is contained in $H^1\big(G, A_s(F_\infty)\big)$. Now since $\Ap^t(F_\infty)=0$ by hypothesis, and $A$ and $A^t$ are isogenous, it follows that $\Ap(F_\infty)$ is finite which in turn implies that $A_s(F_\infty)$ is also finite. Consequently, the kernel of the map $r$ is finite. It follows from this that $R(A_s/F)$ is also finite. Now, combining this latter finiteness observation with a similar argument to that in \cite[Lemma 3.3]{Ha10}, we obtain the conclusion of lemma.
\epf

We now come to the proof of Theorem \ref{no finite submodule theorem}.

\bpf[Proof of Theorem \ref{no finite submodule theorem}]
By Lemma \ref{twist}, we can find a $s\in\Zp$ such that $H^2(G_S(F),A_s)=0$. Now since $A^t(F_\infty)[p^\infty]=0$ by hypothesis, it follows that we also have $A^t_{1-s}(F_\infty) =0$. We can now imitate the proof of Proposition \ref{no finite submodule} to establish that $\Sel(A_s/F_\infty)^{\vee}$ has no non-trivial finite $\Zp\ps{G}$-submodule. It then follows from this that $\Sel(A/F_\infty)^{\vee}$ also has no non-trivial finite $\Zp\ps{G}$-submodule.
\epf

\section{Non-ordinary modular forms} \label{non-ordinary section}

Let $f=\sum a_nq^n$ be a normalised new cuspidal modular eigenform of even weight $k\geq 2$, level $N$ and nebentypus $\epsilon$. We shall always assume that the (odd) prime $p$ does not divide the level $N$ and that the coefficient $a_p$ is not a $p$-adic unit (which is to say $f$ is non-ordinary). Let $\mK_f$ be the number field obtained by adjoining all the Fourier coefficients of $f$ to $\Q$. Throughout, we shall fix a prime $\mathfrak{p}$ of $\mK_f$ above $p$. We then let $V_f$ denote the corresponding two-dimensional $\mK_{f,\mathfrak{p}}$-linear Galois representation attached to $f$ in the sense of Deligne \cite{Del}. Writing $\Op=\Op_{\mK_{f,\mathfrak{p}}}$ for the ring of integers of $\mK_{f,\mathfrak{p}}$, we denote by $T_f$ the $\Gal(\bar{\Q}/\Q)$-stable $\Op$-lattice in $V_f$ as defined in \cite[Section 8.3]{K}. One then sets $A_f = V_f/T_f$. Note that $A_f$ is isomorphic to  $\mK_{f,\mathfrak{p}}/\Op\oplus \mK_{f,\mathfrak{p}}/\Op$ as $\Op$-modules.
Fix an integer $j$ such that $1\leq j\leq k-1$. For any $\Gal(\bar{\Q}/\Q)$-module $M$, we shall write $M(j)$ for $M\ot\chi^j$, where $\chi$ is the $p$-adic cyclotomic character of $\Gal(\bar{\Q}/\Q)$.
All of these fit into the following short exact sequence
\[0\lra T_f(j)\stackrel{\iota}{\lra} V_f(j)\stackrel{\pi}{\lra} A_f(j) \lra 0.\]

Let $F$ be a number field such that $p$ splits completely in $F/\Q$, and let $F_\infty$ be a $\Zp$-extension of $F$ which is totally ramified at every prime of $F$ above $p$. A typical (and important) example of $F_\infty$ is the cyclotomic $\Zp$-extension. Another example is the anti-cyclotomic $\Zp$-extension of an imaginary quadratic field, whose class number is not divisible by $p$.

From now on, $S$ will denote a finite set of primes of $F$ containing those dividing $pN$, the ramified primes of $F/\Q$ and all the infinite primes. Let $L$ be a finite extension of $F$ contained in $F_\infty$. Following Bloch-Kato \cite[(3.7.1), (3.7.2)]{BK}, we define
\[ H^1_{f}\big(L_w, V_f(j)\big)=
\begin{cases} \ker\big(H^1(L_w, V_f(j))\lra H^1(L_w, V\ot\mathbb{B}_{\mathrm{cris}})\big) & \text{\mbox{if} $w$
 divides $p$},\\
 \ker\big(H^1(L_w, V_f(j))\lra H^1(L^{ur}_w, V_f(j))\big) & \text{\mbox{if} $w$ does not divide $p$,}
\end{cases} \]
where $\mathbb{B}_{\mathrm{cris}}$ is the ring of periods of Fontaine \cite{Fon94, Fon04} and $L_w^{ur}$ denotes the maximal unramified extension of $L_w$. Set $H^1_f(L_w, A_f(j)) = \pi_* \big(H^1_f(L_w, V_f(j))\big)$.

If $v$ is a prime of $F$, we then write
\[J_v\big(A_f(j)/L\big) = \bigoplus_{w|v}\frac{H^1(L_w,A_f(j))}{H^1_f(L_w,A_f(j))}\]
and define
\[J_v\big(A_f(j)/F_\infty\big) = \ilim_L J_v(A_f(j)/L),\]
where the limit is taken over finite extension $L$ of $F$ contained in $F_\infty$.

  The Selmer group of $f$ over $F_\infty$ is defined by
\[\Sel\big(A_f(j)/F_\infty\big) = \ker\Big(H^1\big(G_S(F_\infty), A_f(j)\big)\lra \bigoplus_{v\in S}J_v\big(A_f(j)/F_\infty\big)\Big).\]
This has a natural $\Op\ps{G}$-module structure, where $G=\Gal(F_\infty/F)$. Denote by $X\big(\Af/F_\infty\big)$ the Pontryagin dual of $\Sel\big(A_f(j)/F_\infty\big)$.

 In subsequent discussion, we write $T_f^*= \Hom_{\mK_{f,\mathfrak{p}}}(A_f, \mK_{f,\mathfrak{p}}/\Op)$, $V_f^*= \Hom_{\mK_{f,\mathfrak{p}}}(V_f, \mK_{f,\mathfrak{p}})$ and $A_f^*= \Hom_{\mK_{f,\mathfrak{p}}}(T_f, \mK_{f,\mathfrak{p}}/\Op)$.
We can now prove the following analog of Theorem \ref{WLSch}.

\bt \label{WLSch modular form}
 Retain the settings as above. Suppose further that for every prime $v\in S$ that splits completely in $F_\infty/F$, one has $H^0(F_v, V_f^*(1-j))=0$.
 Then we have $H^2(G_S(F_{\infty}),\Af\big)=0$ if and only if
\[\corank_{\Op\ps{G}}\Big(\Sel(\Af/F_{\infty})\Big) = |F:\Q|.\]
\et

\bpf
The proof is essentially similar to that in Theorem \ref{WLSch}. The only thing which perhaps requires additional attention is to establish the analog of Lemma \ref{WLSch lemma}. Since $p$ split completely in $F/\Q$ and $F_\infty/F$ is totally ramified at all primes of $F$ above $p$, we see that $F_{\infty, w}$ is the cyclotomic $\Zp$-extension of $F_v$, where $v$ is a prime of $F$ above $p$ and $w$ is a prime of $F_\infty$ above $v$. Therefore, we may apply \cite[Theorem 0.6]{PR00JAMS} to conclude that $J_v\big(\Af/F_\infty\big)=0$. For a prime $v$ that does not divide $p$ and does not split completely, the same argument in Lemma \ref{WLSch lemma} shows that $J_v\big(\Af/F_\infty\big)$ is cotorsion over $\Zp\ps{G}$.

It therefore remains to analyse $J_v\big(\Af/F_\infty\big)$ for those primes not above $v$ at which $v$ split completely in $F_\infty/F$. In this case, it follows from the hypothesis $H^0(F_v, V_f^*(1-j))=0$ that we have a $\Op$-torsion module $H^0(F_v,A_f^*(1-j))$ injecting into $H^1(F_v, T_f^*(1-j))$. As the latter is finitely generated over $\Op$, we have that
$H^0(F_v, A_f^*(1-j))$ is finite. It then follows from this observation and \cite[Proposition 4.3]{Kid} that $J_v\big(\Af/F_\infty\big)$ is cotorsion over $\Zp\ps{G}$. (Actually, Kidwell works with ordinary modular form. But since this part of the argument is concerned with primes outside $p$, one can check easily that his argument carries over to this non-ordinary context.)
The remainder of the proof then proceeds as in Theorem \ref{WLSch}.
\epf

The next result is an analog of Proposition \ref{surjective} which establishes the exactness of the defining sequence of the Selmer group of our non-ordinary modular form.

\bp \label{surjective modular form}
 Retain the assumptions of Theorem \ref{WLSch modular form}.
 Assume that $H^2(G_{S}(F_\infty),\Af)=0$ and  $\Af(F_\infty)$ is finite. Then we have a short exact sequence
 \[ 0\lra \Sel(\Af/F_\infty)\lra H^1\big(G_S(F_{\infty}),\Af\big)\lra \bigoplus_{v\in S}J_v(\Af/F_\infty) \lra 0.\]
\ep

\bpf
Again, the proof proceeds similarly to that in Proposition \ref{surjective}.
\epf

\br
The condition of $\Af(F_\infty)$ being finite is rather mild. In fact, this group is even trivial in many cases (see \cite[Lemma 4.4, and Assumptions 1 and 2]{Lei11})
\er

We end with the following analog of Theorem \ref{no finite submodule theorem} which also refines \cite[Theorem 3.5]{LeiSu}.

\bt \label{no finite submodule mf}
 Retain the assumptions of Theorem \ref{WLSch modular form}. Suppose that $H^2(G_{S}(F_\infty),\Af)=0$ and that $\Af(F_\infty)=0$.
 Then $X(\Af/F_\infty)$ has no non-trivial finite $\Op\ps{G}$-submodule.
\et

\bpf
This follows from a similar argument to that in Theorem \ref{no finite submodule theorem}.
\epf

\end{document}